\documentclass[12pt]{amsart}
\usepackage{amsmath,amstext,amsfonts,amsthm,amssymb}
\usepackage{eucal}
\usepackage{graphicx}
\renewcommand{\subsubsection}[1]{\addtocounter{subsubsection}{1}
{\ \\[3pt]\bf \thesubsubsection. \  #1} }

\swapnumbers

{  \theoremstyle{definition}

}
\newcommand{\bea}{\begin{eqnarray*}}
\newcommand{\eea}{\end{eqnarray*}}
\newcommand{\bean}{\begin{eqnarray}}
\newcommand{\eean}{\end{eqnarray}}

\newcommand{\BZ}{\mathbb{Z}}

\begin{document}


\centerline{FORMULE DE FATEEV} 

\bigskip\bigskip

\bigskip\bigskip

\centerline{V\'eronique Cohen-Aptel}

\vspace{3cm}

\centerline{\bf \S 1. Introduction}

\bigskip\bigskip

On propose une preuve directe d'une formule, due \`a V.Fateev et collaborateurs, sur les 
produits des valeurs de la fonction Gamma li\'es aux syst\`emes de racines. 

Soit $R\subset V$ un syst\`eme de racines fini r\'eduit irr\'eductible de rang $r$ dans un espace vectoriel r\'eel $V$ de dimension $r$; on munit $V$ d'un produit scalaire $W$-invariant $(.|.)$, $W$ \'etant le groupe 
de Weyl et on id\'entifie $V$ \`a son dual \`a l'aide de ce produit, donc les racines duales $\alpha^\vee = 
2\alpha/(\alpha|\alpha)$ appartiennent \`a $V$. 

On utilisera les notations standardes de [B]. Choisissons une base $\{\alpha_i\}_{1\leq i \leq r}$ de racines 
simples. Suivant l'usage on dit que $R$ est simpl\'ement lac\'e si la matrice de Cartan 
$A = ((\alpha_i|\alpha_j^\vee))_{i,j}$ est sym\'etrique. 

Soit 
$$
\theta = \sum_{i=1}^r\ n_i\alpha_i
$$
la plus longue racine. On pose $\alpha_0 = - \theta,\ n_0 = 1$. Les nombres $n_i$ co\"\i ncident 
avec les marques de Kac du graphe de Dynkin affine de $R^{(1)}$ dans la Table Aff 1, [K]. 


Soit $h = \sum_{i=0}^r\ n_i$ le nombre de Coxeter de $R$.  

On pose  
$$
\rho = \frac{1}{2}\sum_{\alpha > 0}\ \alpha,\ \rho^\vee = \frac{1}{2}\sum_{\alpha > 0}\ \alpha^\vee
$$
On d\'efinit une fonction m\'eromorphe
$$
\gamma(x) = \frac{\Gamma(x)}{\Gamma(1-x)}
$$
Dans un article [F] une formule remarquable suivante a \'et\'e d\'ecouverte (cf. formule (66)):  

{\bf 1.1. Th\'eor\`eme} (Fateev). {\it Supposons que $R$ soit simpl\'ement lac\'e. Pour $1\leq i \leq r$, posons
$$
\gamma(R,\alpha_i) := \prod_{\alpha > 0}\ \gamma((\alpha|\rho)/h)^{-(\alpha_i|\alpha)},
$$
$$
k(R) := \prod_{i=1}^r\ n_i^{n_i}
$$
Alors pour tout $i$ 
$$
\gamma(R,\alpha_i) = k(R)^{-1/h}n_i
\eqno{(F)}
$$
}

Afin d'\'enoncer la formule pour les syst\`emes de racines pas forcement simpl\'ement lac\'es, d\'efinissons les nombres
$$
\gamma'(R,\alpha_i) := \prod_{\alpha > 0}\ \gamma((\alpha|\rho^\vee)/h)^{-(\alpha_i|\alpha^\vee)},
$$
$$
\gamma''(R,\alpha_i) := \prod_{\alpha > 0}\ \gamma((\alpha|\rho)/h^\vee)^{-(\alpha_i^\vee|\alpha)},
$$
$$
k'(R) := \prod_{i=0}^r\ n_i^{\vee n_i}, 
$$
$$
k''(R) := \prod_{i=0}^r\ (n_i^{\vee\vee})^{n_i^\vee} 
$$
o\`u
$$
n_i^\vee = \frac{(\alpha_i|\alpha_i)n_i}{2},\  
$$
$$
n_i^{\vee\vee} := \frac{(\alpha_i|\alpha_i)n_i^\vee}{2}, 
$$
pour $0\leq i \leq n$, 
$$
h^\vee = \sum_{i=0}^n\ n_i^\vee.
$$

Si $R$ est simpl\'ement lac\'e alors $\gamma'(R,\alpha_i) = \gamma''(R,\alpha_i) 
= \gamma(R,\alpha_i)$ et $k'(R) = k''(R) = k(R)$. 

{\bf 1.2. Th\'eor\`eme}, [ABFKR]. {\it Pour tout $i$, 

(i)
$$
\gamma'(R,\alpha_i) = k'(R)^{-1/h}n_i^\vee
\eqno{(F')}
$$

(ii)
$$
\gamma''(R,\alpha_i) = k''(R)^{-1/h^\vee}n_i^{\vee\vee}
\eqno{(F'')}
$$}

Les physiciens d\'eduisent ces formules des arguments tr\`es inter\'essants 
et tr\`es indirects; ils s'appuyient sur l'Ansatz de Bethe. Dans cette note 
on les v\'erifie n'en utilisant que deux propri\'etes fondamentales 
de la fonction Gamma:
$$
\Gamma(x)\Gamma(1-x) = \frac{\pi}{\sin(\pi x)}
\eqno{(C)}
$$
et
$$
\prod_{i=0}^{n-1} \Gamma(x + i/n) = (2\pi)^{(n-1)/2}n^{-nx + 1/2}\Gamma(nx)
\eqno{(M)}
$$
Le fait que ces formules sont d\'eduisables \`a partir de relations $(C)$ et $(M)$ donne lieu 
aux cons\'equences arithm\'etiques, li\'ees aux sommes de Jacobi. Pour leur discussion le lecteur 
est renvoy\'e \`a [CS].

\bigskip\bigskip

\centerline{\bf \S 2. Preuves: cas simpl\'ement lac\'e}

\bigskip\bigskip

On va v\'erifier la formule de Fateev cas par cas. On pr\'esente les nombres 
$\gamma(R,\alpha_i)$ sous une forme explicite
$$
\gamma(R,\alpha_i) = \prod_{j=1}^{h-1}\ \gamma(j/h)^{a_j} = 
\prod_{j=1}^{h-1}\ \Gamma(j/h)^{b_j},
$$
en employant les Planches de [B]. Ensuite, \`a l'aide de (C) et (M) on montre qu'ils sont \'egaux \`a 
$n_ik(R)^{-1/h}$. 

On notera la formule (M) pour $n$ specifique par $M(n)$. 

On des formules \'evidentes
$$
\gamma(x)\gamma(1 - x) = 1;\ \gamma(1/2) = 1
\eqno{(2.1)}
$$
Pour simplifier (et visualiser) l'\'ecriture, on utilisera la notation abr\'eg\'ee:
$$
\{a\}: = \gamma(a/h),\ [a] = \Gamma(a/h), 
\eqno{(2.2)}
$$
o\`u pour chaque s\'erie $R$ ci-dessous $h = h(R)$. 

\bigskip

{\it Syst\`emes de type $A_n$, $n\geq 2$}

\bigskip

$h = n+1$. $\theta = \sum_{i=1}^n\ \alpha_i$, donc tout $n_i = 1$ et $k(A_n) = 1$. 

D'une autre part, 
$$
\gamma(A_n,\alpha_i) = \gamma(i/(n+1))^{-1}\gamma((n+1-i)/(n+1))^{-1} = 1,
$$
d'o\`u le r\'esultat. 

\bigskip

{\it Syst\`emes de type $D_n$, $n\geq 3$}

\bigskip

$h = 2n - 2$. 
$$
\theta = \alpha_1 + 2\alpha_2 + \ldots + 2\alpha_{n-2} + \alpha_{n-1} 
+ \alpha_{n},
$$
d'o\`u
$$
k(D_n) = 2^{2n-6},\ k(D_n)^{-1/h} = 2^{-(2n-6)/(2n-2)}
$$
$$
\gamma(D_n,\alpha_1) = \frac{\gamma((n-2)/(2n-2))}
{\gamma(1/(2n-2))\gamma((n-1)/(2n-2))\gamma((2n-4)/(2n-2))} = 
$$
$$
= \frac{\{n-2\}}{\{1\}\{2n-4\}} = \frac{[2n-3][n-2][2]}{[1][n][2n-4]}
$$
D'apr\`es (D)
$$
[n-2][2n-3] = \pi^{1/2}2^{2/(2n-2)}[2n-4]
$$
et
$$
[1][n] = \pi^{1/2}2^{(2n-4)/(2n-2)}[2],
$$
d'o\`u
$$
\gamma(D_n,\alpha_1) = 2^{(6-2n)/(2n-2)} = n_1k(D_n)^{-1/h}
$$
Pour $2\leq i\leq n-2$, 
$$
\gamma(D_n,\alpha_i) = \frac{\{n-i-1\}\{2n-2i\}}
{\{i\}\{n-i\}\{2n-2i-2\}\{2n-i-1\}} =
$$
$$
= \frac{[2n-i-2][n+i-2][2i][n-i-1][2n-2i][i-1]}
{[i][n-i][2n-2i-2][n+i-1][2i-2][2n-i-1]} = 
$$
$$
= 2^{4/(2n-2)} = n_ik(D_n)^{-1/h}
$$
(on applique (D) avec $x = i/h, (i-1)/h, (n-i)/h, (n-i-1)/h$). Finalement, pour $i = n-1, n$, 
$$
\gamma(D_n,\alpha_i) = \frac{\{2\}}
{\{1\}\{n-1\}\{n\}} = \gamma(D_n,\alpha_1)
$$
$$
= 2^{(6-2n)/(2n-2)} = n_ik(D_n)^{-1/h}.
$$
Ceci \'etablit (F) pour $D_n$. 

\bigskip

{\it Syst\`eme de type $E_6$}

\bigskip

$h = 12$. 
$$
\theta = \alpha_1 + 2\alpha_2 + 2\alpha_3 + 3\alpha_4 + 2\alpha_5 + \alpha_6,
$$
d'o\`u
$$
k(E_6) = 2^6 3^3
$$

$$
\gamma(E_{6},\alpha_{1})=
\biggl\{\gamma(\frac{1}{12})\gamma(\frac{8}{12})\gamma(\frac{3}{12})^{-1}\biggr\}^{-1} =     \biggl\{\gamma(\frac{1}{12})\gamma(\frac{2}{3})\gamma(\frac{1}{4})^{-1}\biggr\}^{-1} = 
$$
$$
=\biggl\{\frac{\Gamma(\frac{1}{12})}{\Gamma(\frac{11}{12})}\frac{\Gamma(\frac{8}{12})}
{\Gamma(\frac{4}{12})}\frac{\Gamma(\frac{9}{12})}{\Gamma(\frac{3}{12})}\biggr\}^{-1}
$$
or d'apr\`es la formule du produit  sur la fonction $\Gamma$ avec $n=3$,
$$ 
2\pi 3^{1/4}\Gamma(\frac{3}{12})=\Gamma(\frac{1}{12})\Gamma(\frac{5}{12})\Gamma(\frac{9}{12})
$$
Donc :
$$
\gamma(E_{6}, \alpha_{1})^{-1}=\frac{2\pi 3^{1/4}\Gamma(\frac{8}{12})}{\Gamma(\frac{5}{12})\Gamma(\frac{11}{12})\Gamma(\frac{4}{12})}
$$
or d'apr\`es la formule de duplication : 
$$
2^{1/6}\pi^{1/2}\Gamma(\frac{5}{6})=\Gamma(\frac{5}{12})\Gamma(\frac{11}{12})
$$
et 
$$
2^{1/3}\pi^{1/2}\Gamma(\frac{2}{3})=\Gamma(\frac{1}{3})\Gamma(\frac{5}{6})
$$

Donc :
$$
\gamma(E_{6}, \alpha_{1})^{-1}=\frac{2\pi 3^{1/4}\Gamma(\frac{2}{3})}{2^{1/6}\pi^{1/2}\Gamma(\frac{5}{6})\Gamma(\frac{1}{3})}=2^{1/2}3^{1/4} 
= n_1^{-1}k^{1/h}
$$

$$
\gamma(E_{6}, \alpha_{2})^{-1}
=\gamma(\frac{5}{12})\gamma(\frac{1}{3})\gamma(\frac{1}{4})^{-1}
$$
or : 
$$ 
2\pi 3^{1/4}\Gamma(\frac{3}{12})=\Gamma(\frac{1}{12})\Gamma(\frac{5}{12})\Gamma(\frac{9}{12}), 
$$
$$ 2^{5/6}\pi ^{1/2}\Gamma(\frac{1}{6})=\Gamma(\frac{1}{12})\Gamma(\frac{7}{12}), 
$$
$$
\Gamma(\frac{1}{6})\Gamma(\frac{5}{6})=2\pi
$$
et
$$2^{1/3}\pi^{1/2}\Gamma(\frac{2}{3})=\Gamma(\frac{1}{3})\Gamma(\frac{5}{6})
$$
Donc : 
 
$$
\gamma(E_{6}, \alpha_{2})^{-1}=
\frac{ 2\pi 3^{1/4}\Gamma(\frac{1}{3})}{\Gamma(\frac{1}{12})\Gamma(\frac{7}{12})\Gamma(\frac{2}{3})}
=\frac{ 2\pi 3^{1/4}\Gamma(\frac{1}{3})}{2^{5/6}\pi ^{1/2}\Gamma(\frac{1}{6})\Gamma(\frac{2}{3})} = 
$$
$$
=\frac{ 2\pi 3^{1/4}\Gamma(\frac{5}{6})\Gamma(\frac{1}{3})}{2^{5/6}\pi ^{1/2}2\pi \Gamma(\frac{2}{3})}
=2^{-1/2}3^{1/4} = n_2^{-1}k^{1/h}
$$

$$
\gamma(E_{6}, \alpha_{4})^{-1}= 
\gamma(\frac{1}{4})^{3}\gamma(\frac{1}{3})^{-1}\gamma(\frac{1}{12})^{-1}\gamma(\frac{5}{12})^{-2}
$$
or :
$$ 
2\pi 3^{-1/4}\Gamma(\frac{3}{4})=\Gamma(\frac{1}{4})\Gamma(\frac{7}{12})\Gamma(\frac{11}{12}), 
$$
$$ 
2^{1/6}\pi ^{1/2}\Gamma(\frac{10}{12})=\Gamma(\frac{5}{12})\Gamma(\frac{11}{12})
$$
et
$$
2^{1/3}\pi^{1/2}\Gamma(\frac{2}{3})=\Gamma(\frac{1}{3})\Gamma(\frac{5}{6})
$$
Donc: 
$$
\gamma(E_{6}, \alpha_{4})^{-1}= 2^{1/2}3^{-3/4} = n_4^{-1}k^{1/h}
$$
Ensuite, on a 
$$
\gamma(E_{6}, \alpha_{3}) = \gamma(E_{6}, \alpha_{5}) = \gamma(E_{6}, \alpha_{2}) 
$$
et
$$
\gamma(E_{6}, \alpha_{6}) = \gamma(E_{6}, \alpha_{1}), 
$$
ce qui implique (F) pour $E_6$. 

\bigskip

{\it Syst\`eme de type $E_7$}

\bigskip

$h = 18$. 
$$
\theta = 2\alpha_1 + 2\alpha_2 + 3\alpha_3 + 4\alpha_4 + 3\alpha_5 + 2\alpha_6 + 
\alpha_7,
$$
d'o\`u
$$
k(E_7) = 2^6 3^2 4 = 2^8 3^2
$$

$$
\gamma(E_{7}, \alpha_{1})^{-1}=\gamma(\frac{6}{18})\gamma(\frac{8}{18})\gamma(\frac{10}{18})
\gamma(\frac{16}{18})^{-1}\gamma(\frac{17}{18})\gamma(\frac{3}{18})^{-1}\gamma(\frac{5}{18})^{-1}\gamma(\frac{1}{18}) = 
$$
$$
= \gamma(\frac{2}{18})\gamma(\frac{3}{18})^{-1}\gamma(\frac{5}{18})^{-1}\gamma(\frac{6}{18}) = 
$$
$$
= \gamma(\frac{1}{9})\gamma(\frac{1}{6})^{-1}\gamma(\frac{5}{18})^{-1}\gamma(\frac{1}{3})
$$

Or utilisant : 
$$\Gamma(\frac{6}{18})= (2\pi)^{-1}3^{-1/6} {\Gamma(\frac{2}{18})\Gamma(\frac{8}{18})\Gamma(\frac{14}{18})}$$
$$\Gamma(\frac{12}{18})= (2\pi)^{-1}3^{1/6} {\Gamma(\frac{4}{18})\Gamma(\frac{10}{18})\Gamma(\frac{16}{18})}$$
$$\Gamma(\frac{3}{18})= (2\pi)^{-1}3^{-1/3} {\Gamma(\frac{1}{18})\Gamma(\frac{7}{18})\Gamma(\frac{13}{18})}$$
$$\Gamma(\frac{15}{18})= (2\pi)^{-1}3^{1/3} {\Gamma(\frac{5}{18})\Gamma(\frac{11}{18})\Gamma(\frac{17}{18})}$$
et :
$$\Gamma(\frac{2}{18})= \pi^{-1/2}2^{-8/9} \Gamma(\frac{1}{18})\Gamma(\frac{10}{18})$$
$$\Gamma(\frac{16}{18})= \pi^{-1/2}2^{-1/9} \Gamma(\frac{8}{18})\Gamma(\frac{17}{18})$$
$$\Gamma(\frac{4}{18})= \pi^{-1/2}2^{-7/9} \Gamma(\frac{2}{18})\Gamma(\frac{11}{18})$$
$$\Gamma(\frac{14}{18})= \pi^{-1/2}2^{-2/9} \Gamma(\frac{7}{18})\Gamma(\frac{16}{18})$$
On obtient :
$$
\gamma(E_{7}, \alpha_{1})=3^{-1/3}2^{2/9} = n_1k(E_7)^{-1/18}
$$

$$
\gamma(E_{7}, \alpha_{2})^{-1}=\gamma(\frac{1}{18})\gamma(\frac{10}{18})^{-1}
\gamma(\frac{2}{18})^{-1}\gamma(\frac{3}{18})^{-1}\gamma(\frac{7}{18})\gamma(\frac{14}{18})\gamma(\frac{6}{18})\gamma(\frac{5}{18})
$$
 Or utilisant : 
$$\Gamma(\frac{6}{18})= (2\pi)^{-1}3^{-1/6} \Gamma(\frac{2}{18})\Gamma(\frac{8}{18})\Gamma(\frac{14}{18})$$
$$\Gamma(\frac{12}{18})= (2\pi)^{-1}3^{1/6} \Gamma(\frac{4}{18})\Gamma(\frac{10}{18})\Gamma(\frac{16}{18})$$
$$\Gamma(\frac{3}{18})= (2\pi)^{-1}3^{-1/3} \Gamma(\frac{1}{18})\Gamma(\frac{7}{18})\Gamma(\frac{13}{18})$$
$$\Gamma(\frac{15}{18})= (2\pi)^{-1}3^{1/3} \Gamma(\frac{5}{18})\Gamma(\frac{11}{18})\Gamma(\frac{17}{18})$$
et :
$$\Gamma(\frac{2}{18})= \pi^{-1/2}2^{-8/9} \Gamma(\frac{1}{18})\Gamma(\frac{10}{18})$$
$$\Gamma(\frac{16}{18})= \pi^{-1/2}2^{-1/9} \Gamma(\frac{8}{18})\Gamma(\frac{17}{18})$$
$$\Gamma(\frac{4}{18})= \pi^{-1/2}2^{-7/9} \Gamma(\frac{2}{18})\Gamma(\frac{11}{18})$$
$$\Gamma(\frac{14}{18})= \pi^{-1/2}2^{-2/9} \Gamma(\frac{7}{18})\Gamma(\frac{16}{18})$$
$$\Gamma(\frac{8}{18})=\pi^{-1/2}2^{-5/9}\Gamma(\frac{4}{18})\Gamma(\frac{13}{18})$$
$$\Gamma(\frac{10}{18})=\pi^{-1/2}2^{-4/9}\Gamma(\frac{5}{18})\Gamma(\frac{14}{18}),$$
on obtient:  
$$
\gamma(E_{7}, \alpha_{2})=2^{-2/9}3^{1/3}=\gamma(E_{7}, \alpha_{1}) = n_2k(E_7)^{-1/18}
$$

$$
\gamma(E_{7}, \alpha_{3})^{-1}= \gamma(\frac{5}{18})\gamma(\frac{11}{18})\gamma(\frac{16}{18})\gamma(\frac{15}{18})^{-1}\gamma(\frac{8}{18})^{-1}
$$

Or:
$$ 2\pi 3^{-1/3}\Gamma(\frac{15}{18})=\Gamma(\frac{5}{18})\Gamma(\frac{11}{18})\Gamma(\frac{17}{18})$$
$$ 2\pi 3^{1/3}\Gamma(\frac{3}{18})=\Gamma(\frac{1}{18})\Gamma(\frac{7}{18})\Gamma(\frac{13}{18})$$
$$ 2^{8/9}\pi ^{1/2}\Gamma(\frac{1}{9})=\Gamma(\frac{1}{18})\Gamma(\frac{5}{9})$$
$$2^{1/9}\pi^{1/2}\Gamma(\frac{8}{9})=\Gamma(\frac{4}{9})\Gamma(\frac{17}{18})$$

Donc : 
$$
\gamma(E_{7}, \alpha_{3})=2^{-7/9}3^{2/3} = n_3k(E_7)^{-1/18}
$$

$$
\gamma(E_{7}, \alpha_{4})^{-1}=\gamma(\frac{8}{18})\gamma(\frac{12}{18})\gamma(\frac{15}{18})\gamma(\frac{14}{18})^{-1}
\gamma(\frac{1}{18})^{-1}\gamma(\frac{7}{18})\gamma(\frac{3}{18})\gamma(\frac{2}{18})\gamma(\frac{5}{18})^{-1}=
$$
$$
= \gamma(\frac{8}{18})\gamma(\frac{12}{18})\gamma(\frac{14}{18})^{-1}\gamma(\frac{1}{18})^{-1}
\gamma(\frac{7}{18})\gamma(\frac{2}{18})\gamma(\frac{5}{18})^{-1}
$$

Or utilisant : 
$$\Gamma(\frac{6}{18})= (2\pi)^{-1}3^{-1/6} \Gamma(\frac{2}{18})\Gamma(\frac{8}{18})\Gamma(\frac{14}{18})$$
$$\Gamma(\frac{12}{18})= (2\pi)^{-1}3^{1/6} \Gamma(\frac{4}{18})\Gamma(\frac{10}{18})\Gamma(\frac{16}{18})$$
 et :
 $$\Gamma(\frac{2}{18})= \pi^{-1/2}2^{-8/9} \Gamma(\frac{1}{18})\Gamma(\frac{10}{18})$$
$$\Gamma(\frac{16}{18})= \pi^{-1/2}2^{-1/9} \Gamma(\frac{8}{18})\Gamma(\frac{17}{18})$$
$$\Gamma(\frac{4}{18})= \pi^{-1/2}2^{-7/9} \Gamma(\frac{2}{18})\Gamma(\frac{11}{18})$$
$$\Gamma(\frac{14}{18})= \pi^{-1/2}2^{-2/9} \Gamma(\frac{7}{18})\Gamma(\frac{16}{18})$$
$$\Gamma(\frac{8}{18})=\pi^{-1/2}2^{-5/9}\Gamma(\frac{4}{18})\Gamma(\frac{13}{18})$$
$$\Gamma(\frac{10}{18})=\pi^{-1/2}2^{-4/9}\Gamma(\frac{5}{18})\Gamma(\frac{14}{18})$$
 
On obtient  : 
$$
\gamma(E_{7}, \alpha_{4})= 2^{11/9}3^{-1/3} = n_4k(E_7)^{-1/18}
$$

$$
\gamma(E_{7}, \alpha_{7})^{-1}=\gamma(\frac{12}{18})\gamma(\frac{1}{18})\gamma(\frac{4}{18})^{-1} 
$$
En utilisant:
$$\Gamma(\frac{6}{18})= (2\pi)^{-1}3^{-1/6} \Gamma(\frac{2}{18})\Gamma(\frac{8}{18})\Gamma(\frac{14}{18})$$
$$\Gamma(\frac{12}{18})= (2\pi)^{-1}3^{1/6} \Gamma(\frac{4}{18})\Gamma(\frac{10}{18})\Gamma(\frac{16}{18})$$
$$\Gamma(\frac{2}{18})= \pi^{-1/2}2^{-8/9} \Gamma(\frac{1}{18})\Gamma(\frac{10}{18})$$
$$\Gamma(\frac{16}{18})= \pi^{-1/2}2^{-1/9} \Gamma(\frac{8}{18})\Gamma(\frac{17}{18}),$$

on obtient
$$
\gamma(E_{7}, \alpha_{7}) = 2^{-7/9}3^{-1/3} = n_7k(E_7)^{-1/18}
$$

Finalement, 
$$
\gamma(E_{7}, \alpha_{5}) = \gamma(E_{7}, \alpha_{3})
$$
et
$$
\gamma(E_{7}, \alpha_{6}) = \gamma(E_{7}, \alpha_{1}), 
$$
ce qui implique (F) pour $E_7$. 

\bigskip

{\it Syst\`eme de type $E_8$}

\bigskip

$h = 30$. 
$$
\theta = 2\alpha_1 + 3\alpha_2 + 4\alpha_3 + 6\alpha_4 + 5\alpha_5 + 4\alpha_6 + 
3\alpha_7 + 2\alpha_8,
$$
d'o\`u
$$
k(E_8) = 2^2 3^2 4^2\cdot 5\cdot 6 = 2^7 3^3 5
$$

$$
\gamma(E_{8},\alpha_{1})^{-1}=\gamma(\frac{1}{30})\gamma(\frac{23}{30})\gamma(\frac{3}{30})^{-1}\gamma(\frac{5}{30})^{-1}
\gamma(\frac{16}{30})^{-1}\gamma(\frac{8}{30})\gamma(\frac{12}{30})\gamma(\frac{10}{30})
$$

Or utilisant : " la table de 5 "
$$\Gamma(\frac{5}{30})= (2\pi)^{-2}5^{-1/3} \Gamma(\frac{1}{30})\Gamma(\frac{7}{30})\Gamma(\frac{13}{30})\Gamma(\frac{19}{30})\Gamma(\frac{25}{30})$$
$$\Gamma(\frac{10}{30})= (2\pi)^{-2}5^{-1/6} \Gamma(\frac{2}{30})\Gamma(\frac{8}{30})\Gamma(\frac{14}{30})\Gamma(\frac{20}{30})\Gamma(\frac{26}{30})$$
$$\Gamma(\frac{20}{30})= (2\pi)^{-2}5^{1/6} \Gamma(\frac{4}{30})\Gamma(\frac{10}{30})\Gamma(\frac{16}{30})\Gamma(\frac{22}{30})\Gamma(\frac{28}{30})$$
$$\Gamma(\frac{25}{30})= (2\pi)^{-2}5^{1/3} \Gamma(\frac{5}{30})\Gamma(\frac{11}{30})\Gamma(\frac{17}{30})\Gamma(\frac{23}{30})\Gamma(\frac{29}{30})$$
"la table de 3 : "
$$\Gamma(\frac{3}{30})= (2\pi)^{-1}3^{-2/5} \Gamma(\frac{1}{30})\Gamma(\frac{11}{30})\Gamma(\frac{21}{30})$$
$$\Gamma(\frac{6}{30})= (2\pi)^{-1}3^{-3/10} \Gamma(\frac{2}{30})\Gamma(\frac{12}{30})\Gamma(\frac{22}{30})$$
$$\Gamma(\frac{9}{30})= (2\pi)^{-1}3^{-1/5} \Gamma(\frac{3}{30})\Gamma(\frac{13}{30})\Gamma(\frac{23}{30})$$
$$\Gamma(\frac{12}{30})= (2\pi)^{-1}3^{-1/10} \Gamma(\frac{4}{30})\Gamma(\frac{14}{30})\Gamma(\frac{24}{30})$$
$$\Gamma(\frac{18}{30})= (2\pi)^{-1}3^{1/10} \Gamma(\frac{6}{30})\Gamma(\frac{16}{30})\Gamma(\frac{26}{30})$$
$$\Gamma(\frac{24}{30})= (2\pi)^{-1}3^{1/10} \Gamma(\frac{8}{30})\Gamma(\frac{18}{30})\Gamma(\frac{28}{30})$$
$$\Gamma(\frac{27}{30})= (2\pi)^{-1}3^{2/5} \Gamma(\frac{9}{30})\Gamma(\frac{19}{30})\Gamma(\frac{29}{30})$$
 et : " la table de 2 "
$$\Gamma(\frac{2}{30})= \pi^{-1/2}2^{-14/15} \Gamma(\frac{1}{30})\Gamma(\frac{16}{30})$$
$$\Gamma(\frac{4}{30})= \pi^{-1/2}2^{-13/15} \Gamma(\frac{2}{30})\Gamma(\frac{17}{30})$$
$$\Gamma(\frac{8}{30})= \pi^{-1/2}2^{-11/15} \Gamma(\frac{4}{30})\Gamma(\frac{19}{30})$$
$$\Gamma(\frac{14}{30})= \pi^{-1/2}2^{-8/15} \Gamma(\frac{7}{30})\Gamma(\frac{22}{30})$$
$$\Gamma(\frac{16}{30})= \pi^{-1/2}2^{-7/15} \Gamma(\frac{8}{30})\Gamma(\frac{23}{30})$$
$$\Gamma(\frac{22}{30})= \pi^{-1/2}2^{-4/15} \Gamma(\frac{11}{30})\Gamma(\frac{26}{30})$$
$$\Gamma(\frac{26}{30})= \pi^{-1/2}2^{-2/15} \Gamma(\frac{13}{30})\Gamma(\frac{28}{30})$$
$$\Gamma(\frac{28}{30})= \pi^{-1/2}2^{-1/15} \Gamma(\frac{14}{30})\Gamma(\frac{29}{30})$$
et donc

$$\frac{\Gamma(\frac{10}{30})}{\Gamma(\frac{20}{30})}=\Gamma(\frac{2}{30})\Gamma(\frac{8}{30})\Gamma(\frac{14}{30})\Gamma(\frac{26}{30})(2\pi)^{-2}5^{-1/6}.$$
$$\frac{\Gamma(\frac{25}{30})}{\Gamma(\frac{5}{30})}=(2\pi)^{-2}5^{1/3}\Gamma(\frac{11}{30})\Gamma(\frac{17}{30})\Gamma(\frac{23}{30})\Gamma(\frac{29}{30}).$$
$$\frac{\Gamma(\frac{27}{30})}{\Gamma(\frac{3}{30})}=\Gamma(\frac{13}{30})\Gamma(\frac{23}{30})\Gamma(\frac{19}{30})\Gamma(\frac{29}{30})(2\pi)^{-2}3^{1/5}.$$
$$\frac{\Gamma(\frac{12}{30})}{\Gamma(\frac{18}{30})}= \Gamma(\frac{4}{30})\Gamma(\frac{14}{30})\Gamma(\frac{8}{30})\Gamma(\frac{28}{30})(2\pi)^{-2}3^{1/5}.$$

On obtient:
$$
\gamma(E_{8},\alpha_{1})^{-1}=3^{2/5}5^{1/6}(2^{-4}\pi^{-4})^{2}\Gamma(\frac{1}{30})
\Gamma(\frac{23}{30})^{3}\Gamma(\frac{14}{30})^{3}\Gamma(\frac{8}{30})^{3}\times
$$
$$\times
\Gamma(\frac{2}{30})\Gamma(\frac{26}{30})\Gamma(\frac{11}{30})\Gamma(\frac{17}{30})\Gamma(\frac{13}{30})\Gamma(\frac{19}{30})\Gamma(\frac{29}{30})
\times
$$
$$
\times\Gamma(\frac{4}{30})\Gamma(\frac{28}{30})\Gamma(\frac{7}{30})^{-1}\Gamma(\frac{16}{30})^{-1}\Gamma(\frac{22}{30})^{-1}
$$

Puis on remplace, gr\^ace \`a la "table de 2",  chaque expression comme : $$\Gamma(\frac{1}{30}), \Gamma(\frac{23}{30}), \Gamma(\frac{11}{30}), \Gamma(\frac{17}{30}),\Gamma(\frac{13}{30}), \Gamma(\frac{19}{30}), \Gamma(\frac{29}{30}), \Gamma(\frac{7}{30})$$
On a: 
$$\gamma(E_{8},\alpha_{1})^{-1} = 2^{-62/15}3^{2/5}5^{1/6}\pi^{-4}\Gamma(\frac{2}{30})\Gamma(\frac{4}{30})\Gamma(\frac{8}{30})\Gamma(\frac{14}{30})\Gamma(\frac{16}{30})\Gamma(\frac{22}{30})\Gamma(\frac{26}{30})\Gamma(\frac{28}{30})$$
sachant que :  gr\`ace \`a la "table de 5" 
$$
\Gamma(\frac{10}{30})×\Gamma(\frac{20}{30})= (2\pi)^{-4} \Gamma(\frac{2}{30})\Gamma(\frac{8}{30})\Gamma(\frac{14}{30})\Gamma(\frac{20}{30})\Gamma(\frac{26}{30})\Gamma(\frac{4}{30})\Gamma(\frac{10}{30})\Gamma(\frac{16}{30})\Gamma(\frac{22}{30})\Gamma(\frac{28}{30})
$$ 
on a alors que: 
$$
\Gamma(\frac{2}{30})\Gamma(\frac{4}{30})\Gamma(\frac{8}{30})\Gamma(\frac{14}{30})\Gamma(\frac{16}{30})\Gamma(\frac{22}{30})\Gamma(\frac{26}{30})\Gamma(\frac{28}{30}) =( 2\pi)^{4}
$$  
D'o\`u: 
$$
\gamma(E_{8}, \alpha_{1})=2^{2/15}3^{-2/5}5^{-1/6} = n_1k(E_8)^{-1/30}
$$

******************

$$
\gamma(E_{8},\alpha_{2})^{-1}=\gamma(\frac{1}{30})\gamma(\frac{6}{30})\gamma(\frac{7}{30})\gamma(\frac{8}{30})\gamma(\frac{17}{30})\gamma(\frac{15}{30})\gamma(\frac{24}{30})\times
$$
$$
\times \gamma(\frac{2}{30})^{-1}\gamma(\frac{3}{30})^{-1}\gamma(\frac{10}{30})^{-1}\gamma(\frac{12}{30})^{-1}\gamma(\frac{21}{30})^{-1}
$$

Or r\'eutilisant   "la table de 3 et de 5", on a :
$$\frac{\Gamma(\frac{20}{30})}{\Gamma(\frac{10}{30})}= (2\pi)^{-2}5^{1/6} \Gamma(\frac{4}{30})\Gamma(\frac{16}{30})\Gamma(\frac{22}{30})\Gamma(\frac{28}{30})$$
$$\frac{\Gamma(\frac{27}{30})}{\Gamma(\frac{3}{30})}=\Gamma(\frac{13}{30})\Gamma(\frac{23}{30})\Gamma(\frac{19}{30})\Gamma(\frac{29}{30})(2\pi)^{-2}3^{1/5}.$$
$$\frac{\Gamma(\frac{9}{30})}{\Gamma(\frac{21}{30})}=(\Gamma(\frac{7}{30})\Gamma(\frac{17}{30})\Gamma(\frac{19}{30})\Gamma(\frac{29}{30}))^{-1}(2\pi)^{2}3^{-3/5}.$$
$$\frac{\Gamma(\frac{18}{30})}{\Gamma(\frac{12}{30})}=( \Gamma(\frac{4}{30})\Gamma(\frac{14}{30})\Gamma(\frac{8}{30})\Gamma(\frac{28}{30}))^{-1}(2\pi)^{2}3^{-1/5}.$$

On obtient:
$$
\gamma(E_{8}, \alpha_{2})^{-1}= 3^{-3/5}5^{1/6}\frac{\Gamma(\frac{16}{30})\Gamma(\frac{22}{30})\Gamma(\frac{1}{30})\Gamma(\frac{28}{30})}{\Gamma(\frac{14}{30})\Gamma(\frac{29}{30})\Gamma(\frac{22}{30})\Gamma(\frac{2}{30})}
$$
or, d'apr\`es " la table de 2 ":
$$\Gamma(\frac{2}{30})= \pi^{-1/2}2^{-14/15} \Gamma(\frac{1}{30})\Gamma(\frac{16}{30})$$
$$\Gamma(\frac{28}{30})= \pi^{-1/2}2^{-1/15} \Gamma(\frac{14}{30})\Gamma(\frac{29}{30})$$

Il s'en suit: 
$$
\gamma(E_{8}, \alpha_{2})=2^{-13/15}3^{3/5}5^{-1/6} = n_2k(E_8)^{-1/30}
$$

*************

$$
\gamma(E_{8},\alpha_{3})^{-1}
=\gamma(\frac{20}{30})\gamma(\frac{24}{30})\gamma(\frac{13}{30})\gamma(\frac{7}{30})\gamma(\frac{11}{30})\gamma(\frac{15}{30})^{-1}\gamma(\frac{8}{30})^{-1}\gamma(\frac{22}{30})^{-1}
$$

Or utilisant : " la table de 5 et 3 ": 
$$\frac{\Gamma(\frac{20}{30})}{\Gamma(\frac{10}{30})}= (2\pi)^{-2}5^{1/6} \Gamma(\frac{4}{30})\Gamma(\frac{16}{30})\Gamma(\frac{22}{30})\Gamma(\frac{28}{30})$$
$$\frac{\Gamma(\frac{24}{30})}{\Gamma(\frac{6}{30})}= (2\pi)^{2}3^{2/5} \Gamma(\frac{2}{30})\Gamma(\frac{4}{30})\Gamma(\frac{14}{30})\Gamma(\frac{22}{30})$$
on obtient :
 
$$
\gamma(E_{8},\alpha_{3})^{-1}
= 3^{2/5}5^{1/6}\frac{\Gamma(\frac{16}{30})\Gamma(\frac{28}{30})\Gamma(\frac{13}{30})\Gamma(\frac{7}{30})
\Gamma(\frac{11}{30})}{\Gamma(\frac{2}{30})\Gamma(\frac{14}{30})\Gamma(\frac{17}{30})\Gamma(\frac{23}{30})\Gamma(\frac{19}{30})}
$$ 
puis on remplace, gr\^ace \`a la "table de 2",  chaque expression comme : $$\Gamma(\frac{23}{30}), \Gamma(\frac{11}{30}), \Gamma(\frac{17}{30}),\Gamma(\frac{13}{30}), \Gamma(\frac{19}{30}), \Gamma(\frac{7}{30}),
$$
et on obtient 
$$
\gamma(E_{8}, \alpha_{3}) = 2^{17/15}3^{-2/5}5^{-1/6} = n_3k(E_8)^{-1/30}
$$

*****************

$$
\gamma(E_{8},\alpha_{4})^{-1} = \gamma(\frac{2}{30})\gamma(\frac{3}{30})\gamma(\frac{25}{30})\gamma(\frac{11}{30})^{-1}
\gamma(\frac{20}{30})^{-1}\gamma(\frac{24}{30})^{-1}\gamma(\frac{1}{30})^{-1}\gamma(\frac{4}{30})^{-1}
$$

Or r\'eutilisant : la "table de 3 et de 5" et notamment ces r\'esultats: 
$$\frac{\Gamma(\frac{10}{30})}{\Gamma(\frac{20}{30})}=\Gamma(\frac{2}{30})\Gamma(\frac{8}{30})\Gamma(\frac{14}{30})\Gamma(\frac{26}{30})(2\pi)^{-2}5^{-1/6}.$$
$$\frac{\Gamma(\frac{25}{30})}{\Gamma(\frac{5}{30})}=(2\pi)^{-2}5^{1/3}\Gamma(\frac{11}{30})\Gamma(\frac{17}{30})\Gamma(\frac{23}{30})\Gamma(\frac{29}{30}).$$
$$\frac{\Gamma(\frac{3}{30})}{\Gamma(\frac{27}{30})}=(\Gamma(\frac{13}{30})\Gamma(\frac{23}{30})\Gamma(\frac{19}{30})\Gamma(\frac{29}{30})^{-1}(2\pi)^{2}3^{-1/5}.$$
$$\frac{\Gamma(\frac{6}{30})}{\Gamma(\frac{24}{30})}= \Gamma(\frac{2}{30})\Gamma(\frac{4}{30})\Gamma(\frac{14}{30})\Gamma(\frac{22}{30})(2\pi)^{-2}3^{-2/5}, 
$$
on obtient:

$$
\gamma(E_{8}, \alpha_{4})^{-1}
=3^{-3/5}5^{1/6}(2\pi)^{-4}\frac{ \Gamma(\frac{2}{30})^{2}\Gamma(\frac{26}{30})^{2}\Gamma(\frac{14}{30})^{2}\Gamma(\frac{8}{30})\Gamma(\frac{22}{30})\Gamma(\frac{29}{30})\Gamma(\frac{17}{30})}{\Gamma(\frac{13}{30})\Gamma(\frac{1}{30})\Gamma(\frac{28}{30})}
$$

Puis on remplace, grâce \`a la " table de 2 ",  chaque expression comme  
$$
\Gamma(\frac{29}{30}), \Gamma(\frac{17}{30}), \Gamma(\frac{13}{30}),\Gamma(\frac{1}{30}),
$$
d'o\`u 
$$
\gamma(E_{8},\alpha_{4})^{-1}
= 3^{-3/5}5^{1/6}(2\pi)^{-4}2^{-2/15}\Gamma(\frac{2}{30})
\Gamma(\frac{4}{30})\Gamma(\frac{8}{30})\Gamma(\frac{14}{30})\times
$$
$$
\times\Gamma(\frac{16}{30})\Gamma(\frac{22}{30})\Gamma(\frac{26}{30})\Gamma(\frac{28}{30})
$$
Ce qui donne : 
$$
\gamma(E_{8}, \alpha_{4})= 2^{2/15}3^{-2/5}5^{-1/6}(2\pi)^{-4}(2\pi)^{4}=2^{2/15}3^{-2/5}5^{-1/6} = n_4k(E_8)^{-1/30}
$$

*************************

$$
\gamma(E_{8},\alpha_{5})^{-1} =\gamma(\frac{11}{30})\gamma(\frac{7}{30})^{-1}\gamma(\frac{8}{30})^{-1}\gamma(\frac{13}{30})^{-1}\gamma(\frac{16}{30})\gamma(\frac{4}{30})\gamma(\frac{6}{30})^{-1}\gamma(\frac{2}{30})^{-1}\gamma(\frac{20}{30})^{-1}\gamma(\frac{25}{30})^{-2}$$

Or utilisant : les r\'esultats cit\'es pr\'ec\'edemment \`a savoir : 
$$\frac{\Gamma(\frac{10}{30})}{\Gamma(\frac{20}{30})}=\Gamma(\frac{2}{30})\Gamma(\frac{8}{30})\Gamma(\frac{14}{30})\Gamma(\frac{26}{30})(2\pi)^{-2}5^{-1/6}.$$
$$\frac{\Gamma(\frac{5}{30})}{\Gamma(\frac{25}{30})}=(2\pi)^{2}5^{-1/3}(\Gamma(\frac{11}{30})\Gamma(\frac{17}{30})\Gamma(\frac{23}{30})\Gamma(\frac{29}{30}))^{-1}.$$
$$\frac{\Gamma(\frac{24}{30})}{\Gamma(\frac{6}{30})}= (\Gamma(\frac{2}{30})\Gamma(\frac{4}{30})\Gamma(\frac{14}{30})\Gamma(\frac{22}{30}))^{-1}(2\pi)^{2}3^{2/5},$$

on obtient:

$$\gamma(E_{8}, \alpha_{5})^{-1} = 3^{2/5}5^{-5/6}(2\pi)^{4}\frac{\Gamma\frac{16}{30})\Gamma(\frac{28}{30})\Gamma(\frac{8}{30})}{\Gamma(\frac{19}{30}\frac{7}{30})\Gamma(\frac{13}{30})\Gamma(\frac{14}{30})\Gamma(\frac{2}{30}\frac{11}{30})\Gamma(\frac{17}{30})\Gamma(\frac{23}{30})\Gamma(\frac{29}{30})^{2}}$$

puis on remplace, gr\`ace \`a la "table de 2",  chaque expression comme : 
$$
\Gamma(\frac{29}{30}), \Gamma(\frac{17}{30}), \Gamma(\frac{13}{30}),\Gamma(\frac{19}{30}),  \Gamma(\frac{7}{30}),  \Gamma(\frac{11}{30}),  \Gamma(\frac{23}{30})
$$
D'o\`u: 
$$
\gamma(E_{8}, \alpha_{5})=3^{-2/5}5^{5/6}2^{-4+47/15}\pi^{4}\pi^{-4} = 2^{-13/15}3^{-2/5}5^{5/6} = n_5k(E_8)^{-1/30}
$$

*********************

$$
\gamma(E_{8},\alpha_{6})^{-1} =\gamma(\frac{27}{30})\gamma(\frac{20}{30})\gamma(\frac{14}{30})\gamma(\frac{13}{30})\gamma(\frac{6}{30})\gamma(\frac{9}{30})^{-1}\gamma(\frac{15}{30})^{-1}\gamma(\frac{26}{30})^{-1}
$$
Or utilisant les r\'esultats cit\'es pr\'ec\'edemment \`a savoir: 
$$\frac{\Gamma(\frac{20}{30})}{\Gamma(\frac{10}{30})}=(\Gamma(\frac{2}{30})\Gamma(\frac{8}{30})\Gamma(\frac{14}{30})\Gamma(\frac{26}{30}))^{-1}(2\pi)^{2}5^{1/6}.$$
$$\frac{\Gamma(\frac{27}{30})}{\Gamma(\frac{3}{30})}=\Gamma(\frac{13}{30})\Gamma(\frac{23}{30})\Gamma(\frac{19}{30})\Gamma(\frac{29}{30})(2\pi)^{-2}3^{1/5}.$$
$$\frac{\Gamma(\frac{6}{30})}{\Gamma(\frac{24}{30})}= \Gamma(\frac{2}{30})\Gamma(\frac{4}{30})\Gamma(\frac{14}{30})\Gamma(\frac{22}{30})(2\pi)^{-2}3^{-2/5}.$$
$$\frac{\Gamma(\frac{21}{30})}{\Gamma(\frac{9}{30})}=\Gamma(\frac{7}{30})\Gamma(\frac{17}{30})\Gamma(\frac{19}{30})\Gamma(\frac{29}{30})(2\pi)^{-2}3^{3/5},
$$
on obtient:

$$
\gamma(E_{8},\alpha_{6})^{-1} = 3^{2/5}5^{1/6}(2\pi)^{-4}\frac{\Gamma(\frac{13}{30})^{2}\Gamma(\frac{23}{30})\Gamma(\frac{19}{30})^{2}\Gamma(\frac{29}{30})^{2}\Gamma(\frac{4}{30})^{2}\Gamma(\frac{14}{30})\Gamma(\frac{22}{30})\Gamma(\frac{7}{30})}{\Gamma(\frac{8}{30})\Gamma(\frac{16}{30})\Gamma(\frac{26}{30})^{2}}
$$
Puis on remplace, grâce \`a la " table de 2 ",  chaque expression comme : 
$$\Gamma(\frac{29}{30}), \Gamma(\frac{13}{30}),\Gamma(\frac{19}{30}),  \Gamma(\frac{7}{30}),\Gamma(\frac{23}{30})$$

D'o\`u :
$$
\gamma(E_{8}, \alpha_{6})=3^{-2/5}5^{-1/6}(2\pi)^{-4}\pi^{4}2^{-43/15}=2^{17/15}3^{-2/5}5^{-1/6}=\gamma(E_{8}, \alpha_{3}) 
= 
$$
$$
= n_6k(E_8)^{-1/30}
$$

******************

$$
\gamma(E_{8}, \alpha_{7})^{-1}  =\gamma(\frac{9}{30})\gamma(\frac{19}{30})\gamma(\frac{4}{30})^{-1}\gamma(\frac{10}{30})^{-1}\gamma(\frac{14}{30})^{-1}\gamma(\frac{18}{30})^{-1}\gamma(\frac{27}{30})^{-1}
$$
Or utilisant les r\'esultats cit\'es pr\'ec\'edemment \`a savoir : 
$$\frac{\Gamma(\frac{20}{30})}{\Gamma(\frac{10}{30})}=(\Gamma(\frac{2}{30})\Gamma(\frac{8}{30})\Gamma(\frac{14}{30})\Gamma(\frac{26}{30}))^{-1}(2\pi)^{2}5^{1/6}.$$
$$\frac{\Gamma(\frac{3}{30})}{\Gamma(\frac{27}{30})}=(\Gamma(\frac{13}{30})\Gamma(\frac{23}{30})\Gamma(\frac{19}{30})\Gamma(\frac{29}{30}))^{-1}(2\pi)^{2}3^{-1/5}.$$
$$\frac{\Gamma(\frac{12}{30})}{\Gamma(\frac{18}{30})}= \Gamma(\frac{4}{30})\Gamma(\frac{14}{30})\Gamma(\frac{8}{30})\Gamma(\frac{28}{30})(2\pi)^{-2}3^{1/5}.$$
$$\frac{\Gamma(\frac{9}{30})}{\Gamma(\frac{21}{30})}=(\Gamma(\frac{7}{30})\Gamma(\frac{17}{30})\Gamma(\frac{19}{30})\Gamma(\frac{29}{30}))^{-1}(2\pi)^{2}3^{-3/5},
$$

on obtient:

$$
\gamma(E_{8},\alpha_{7})^{-1} = 3^{-3/5}5^{1/6}\frac{\Gamma(\frac{16}{30})\Gamma(\frac{28}{30})\Gamma(\frac{1}{30})}
{\Gamma(\frac{29}{30})\Gamma(\frac{2}{30})\Gamma(\frac{14}{30})}
$$
et sachant que:
$$\Gamma(\frac{2}{30})= \pi^{-1/2}2^{-14/15} \Gamma(\frac{1}{30})\Gamma(\frac{16}{30})$$
$$\Gamma(\frac{28}{30})= \pi^{-1/2}2^{-1/15} \Gamma(\frac{14}{30})\Gamma(\frac{29}{30})$$
il vient: 
$$
\gamma(E_{8}, \alpha_{7}) = 2^{-13/15}3^{3/5}5^{-1/6} = \gamma(E_{8}, \alpha_{2}) = n_7k(E_8)^{-1/30}
$$

*****************

$$
\gamma(E_{8},\alpha_{8})^{-1} =\gamma(\frac{14}{30})\gamma(\frac{18}{30})\gamma(\frac{28}{30})^{-1}\gamma(\frac{5}{30})^{-1}\gamma(\frac{10}{30})\gamma(\frac{9}{30})^{-1}
$$

Or utilisant les r\'esultats cit\'es pr\'ec\'edemment, \`a savoir : 
$$\frac{\Gamma(\frac{10}{30})}{\Gamma(\frac{20}{30})}=\Gamma(\frac{2}{30})\Gamma(\frac{8}{30})\Gamma(\frac{14}{30})\Gamma(\frac{26}{30})(2\pi)^{-2}5^{-1/6}.$$
$$\frac{\Gamma(\frac{25}{30})}{\Gamma(\frac{5}{30})}=(2\pi)^{-2}5^{1/3}\Gamma(\frac{11}{30})\Gamma(\frac{17}{30})\Gamma(\frac{23}{30})\Gamma(\frac{29}{30}).$$
$$\frac{\Gamma(\frac{21}{30})}{\Gamma(\frac{9}{30})}=\Gamma(\frac{7}{30})\Gamma(\frac{17}{30})\Gamma(\frac{19}{30})\Gamma(\frac{29}{30})(2\pi)^{-2}3^{3/5}.$$
$$\frac{\Gamma(\frac{18}{30})}{\Gamma(\frac{12}{30})}= (\Gamma(\frac{4}{30})\Gamma(\frac{14}{30})\Gamma(\frac{8}{30})\Gamma(\frac{28}{30}))^{-1}(2\pi)^{2}3^{-1/5},
$$
on obtient:
$$
\gamma(E_{8},\alpha_{8})^{-1} = 3^{2/5}5^{1/6}(2\pi)^{-4}\frac{\Gamma(\frac{2}{30})\Gamma(\frac{26}{30})\Gamma(\frac{11}{30})\Gamma(\frac{17}{30})^{2}\Gamma(\frac{23}{30})\Gamma(\frac{29}{30})^{2}\Gamma(\frac{7}{30})\Gamma(\frac{19}{30})\Gamma(\frac{14}{30})}{\Gamma(\frac{4}{30})\Gamma(\frac{16}{30})\Gamma(\frac{28}{30})}$$
On remplace, gr\`ace \`a la " table de 2 ", les expressions comme 
$$
\Gamma(\frac{11}{30}),\Gamma(\frac{17}{30}), \Gamma(\frac{23}{30}), \Gamma(\frac{29}{30}),\Gamma(\frac{7}{30})\Gamma(\frac{19}{30}) .
$$
Il vient: 
$$
\gamma(E_{8}, \alpha_{8}) = 2^{2/15}3^{-2/5}5^{-1/6}= \gamma(E_{8}, \alpha_{1}) = n_7k(E_8)^{-1/30}
$$
Ceci \'etablit (F) pour le syst\`eme de type $E_8$ et termine la preuve du Th\'eor\`eme 1.1. $\square$

\bigskip\bigskip

\centerline{\bf \S 3. Preuves: cas non-simpl\'ement lac\'e}

\bigskip\bigskip

Dans ce paragraphe on va d\'emontrer Th\'eor\`eme 1.2. On utilise toujours la notation de 
[B], Planches \`a la fin du livre. 

\bigskip

{\it Syst\`emes de type $B_n$, $n\geq 2$}
 
\bigskip

$h = 2n$. 
$$
\theta = \alpha_1 + 2\alpha_2 + 2\alpha_3 + \ldots + 2\alpha_n
$$
Racines simples:
$$
\alpha_0 = -\epsilon_1 - \epsilon_2, \alpha_i = \epsilon_i - \epsilon_{i+1},\ 1\leq i \leq n-1,\ \alpha_n = \epsilon_n
$$
Il s'en suit:
$$
(n_i)_{0\leq i\leq n} = (1,1,2,2,\ldots, 2, 2),
$$
$$
(n_i^\vee)_{0\leq i\leq n} = (1,1,2,2,\ldots, 2, 1),
$$
$$
(n_i^{\vee\vee})_{0\leq i\leq n} = (1,1,2,2,\ldots, 2, 1/2),
$$
$$
h^\vee = 2n - 1
$$
$$
k'(B_n) = 2^{2n-4},\ k''(B_n) = 2^{2n-5}
$$

$$
\gamma'(B_n,\alpha_1) = \gamma\big(\frac{n-1}{2n}\big)\gamma\big(\frac{1}{2n}\big)^{-1}\gamma\big(\frac{2n-2}{2n}\big) 
= \gamma\big(\frac{n-1}{2n}\big)\gamma\big(\frac{1}{2n}\big)^{-1}\gamma\big(\frac{n-1}{n}\big)^{-1}
$$
En sachant que
$$\Gamma\big(\frac{2}{2n}\big)=\Gamma\big(\frac{1}{2n}\big)\Gamma\big(\frac{n+1}{2n}\big)\pi^{-1/2}2^{\frac{2-2n}{2n}}$$
$$
\Gamma\big(\frac{2n-2}{2n}\big)=\Gamma\big(\frac{l-1}{2n}\big)\Gamma\big(\frac{2n-1}{2n}\big)\pi^{-1/2}2^{\frac{-2}{2n}},
$$
on obtient
$$
\gamma'(B_n,\alpha_1) = 
\pi^{-1/2}2^{-\frac{2}{2n}}\pi^{1/2}2^{\frac{-2+2n}{2n}}=2^{\frac{-4+2n}{2n}} 
= n_1^\vee k'(B_n)^{-1/2n}
$$

$$
\gamma'(B_n,\alpha_n) = \gamma'(B_n,\alpha_1) = n_n^\vee k'(B_n)^{-1/2n}
$$

Pour $2 \leq i \leq n-1$

$$
\gamma'(B_n,\alpha_i) =
\gamma\big(\frac{n-i}{2n}\big)\gamma\big(\frac{n-i+1}{2n}\big)^{-1}\gamma\big(\frac{2n-2i}{2n}\big)^{-1}\times
$$
$$
\times \gamma\big(\frac{2n-2i+2}{2n}\big)\gamma\big(\frac{2n-i+1}{2n}\big)^{-1}\gamma\big(\frac{i}{2n}\big)^{-1}
$$
En sachant que 
$$\Gamma\big(\frac{2i-2}{2n}\big)=\Gamma\big(\frac{i-1}{2n}\big)\Gamma\big(\frac{n+i-1}{2n}\big)\pi^{-1/2}2^{\frac{2i-2-2n}{2n}}$$
$$\Gamma\big(\frac{2n-2i}{2n}\big)=\Gamma\big(\frac{n-i}{2n}\big)\Gamma\big(\frac{2n-i}{2n}\big)\pi^{-1/2}2^{\frac{-2i}{2n}}$$
$$\Gamma\big(\frac{2n-2i+2}{2n}\big)=\Gamma\big(\frac{n-i+1}{2n}\big)\Gamma\big(\frac{2n-i+1}{2n}\big)\pi^{-1/2}2^{\frac{2-2i}{2n}}$$
$$
\Gamma\big(\frac{2i}{2n}\big)=\Gamma\big(\frac{i}{2n}\big)\Gamma\big(\frac{i+n}{2n}\big)\pi^{-1/2}2^{\frac{2i-2n}{2n}},
$$
on obtient
$$
\gamma'(B_n,\alpha_i) = 2^{-\frac{2}{n}} = n_i^\vee k'(B_n)^{-1/2n}
$$
Donc $(F')$ est v\'erifi\'ee pour le syst\`eme du type $B_{n}.$

\bigskip

{\it Formule} $(F'')$


$$
\gamma''(B_n,\alpha_1) = 
\gamma\big(\frac{2}{2n-1}\big)\gamma\big(\frac{1}{2n-1}\big)^{-1}\gamma\big(\frac{2n+1}{2(2n-1)}\big)^{-1} 
= \gamma\big(\frac{n-1}{2n}\big)\gamma\big(\frac{1}{2n}\big)^{-1}\gamma\big(\frac{n-1}{n}\big)^{-1}
$$
En sachant que:
$$\Gamma\big(\frac{2}{2n-1}\big)=\Gamma\big(\frac{1}{2n-1}\big)\Gamma\big(\frac{2n+1}{2(2n-1)}\big)\pi^{-1/2}2^{\frac{3-2n}{2n-1}}$$
$$\Gamma\big(\frac{2n-3}{2n-1}\big)=\Gamma\big(\frac{2n-3}{2(2n-1)}\big)\Gamma\big(\frac{2(2n-2)}{2n-1}\big)\pi^{-1/2}2^{\frac{-2}{2n-1}},
$$
on obtient
$$
\gamma\big(\frac{n-1}{2n}\big)\gamma\big(\frac{1}{2n}\big)^{-1}\gamma\big(\frac{n-1}{n}\big)^{-1} 
= \pi^{-1/2}2^{\frac{2}{2n-1}}\pi^{1/2}2^{\frac{3-2n}{2n-1}}=2^{\frac{5-2n}{2n-1}} = n_1^{\vee\vee} k''(B_n)^{-1/(2n-1)}
$$
Donc la formule $(F'')$ est v\'erifi\'ee pour $i=1$

V\'erifions maintenant la formule $(F'')$ pour $i=n$. 

$$
\gamma''(B_n,\alpha_n) = 
\gamma\big(\frac{1}{2(2n-1)}\big)^{-2}\gamma\big(\frac{1}{2n-1}\big)^{2}\gamma\big(\frac{n}{2n-1}\big)^{-2} 
$$
En sachant que:
$$\Gamma\big(\frac{1}{2n-1}\big)=\Gamma\big(\frac{1}{4l-2}\big)\Gamma\big(\frac{n}{2n-1}\big)\pi^{-1/2}2^{\frac{2-2n}{2n-1}}$$
$$\Gamma\big(\frac{2n-2}{2n-1}\big)=\Gamma\big(\frac{4l-3}{2(2n-1)}\big)\Gamma\big(\frac{l-1}{2n-1}\big)\pi^{-1/2}2^{\frac{-1}{2n-1}},
$$
il vient:
$$
\gamma\big(\frac{1}{2(2n-1)}\big)^{-2}\gamma\big(\frac{1}{2n-1}\big)^{2}\gamma\big(\frac{n}{2n-1}\big)^{-2}= 2^{\frac{6-4n}{2n-1}} 
= n_n^{\vee\vee} k''(B_n)^{-1/(2n-1)} 
$$
Donc la formule $(F'')$ est v\'erifi\'ee pour $i=n$. 

Pour $2 \leqslant i\leqslant n-1$, 
$$
\gamma''(B_n,\alpha_i) = \gamma\big(\frac{2n-2i+1}{2(2n-1)}\big)^{-1}\gamma\big(\frac{2n-2i-1}{2(2n-1)}\gamma\big(\frac{2n-2i-1}{2n-1}\big)^{-1}
\times
$$
$$
\times
\gamma\big(\frac{2n-2i+1}{2n-1}\big)\gamma\big(\frac{2n-i}{2n-1}\big)^{-1}\gamma\big(\frac{i}{2n-1}\big)^{-1} =2^{\frac{4}{2n-1}}
$$
En sachant que:
$$\Gamma\big(\frac{2i-2}{2n-1}\big)=\Gamma\big(\frac{i-1}{2n-1}\big)\Gamma\big(\frac{2n+2i-3}{2(2n-1)}\big)\pi^{-1/2}2^{\frac{2i-1-2n}{2n-1}}$$
$$\Gamma\big(\frac{2n-2i+1}{2n-1}\big)=\Gamma\big(\frac{2n-i}{2n-1}\big)\Gamma\big(\frac{2n-2i+1}{2(2n-1)}\big)\pi^{-1/2}2^{\frac{-2i+2}{2n-1}}$$
$$\Gamma\big(\frac{2n-2i-1}{2n-1}\big)=\Gamma\big(\frac{2n-2i-1}{2n-1}\big)\Gamma\big(\frac{2n-i-1}{2n-1}\big)\pi^{-1/2}2^{\frac{-2i}{2n-1}}$$
$$\Gamma\big(\frac{2i}{2n-1}\big)=\Gamma\big(\frac{i}{2n-1}\big)\Gamma\big(\frac{2i+2n-1}{2(2n-1)}\big)\pi^{-1/2}2^{\frac{2i-2n+1}{2n-1}},
$$
on obtient
$$
\gamma''(B_n,\alpha_i) = 2^{\frac{4}{2n-1}} = n_i^{\vee\vee} k''(B_n)^{-1/(2n-1)}
$$
Donc la formule $(F'')$ est v\'erifi\'ee pour $2\leqslant i \leqslant n-1.$

Donc $(F'')$ est v\'erifi\'ee pour le syst\`eme du type $B_{n}.$

\bigskip

{\it Syst\`emes de type $C_n$, $n\geq 2$}
 
\bigskip

$h = 2n$. 
$$
\theta = 2\alpha_1 + 2\alpha_2 + 2\alpha_3 + \ldots + \alpha_n
$$
Racines simples:
$$
\alpha_0 = - 2\epsilon_1, \alpha_i = \epsilon_i - \epsilon_{i+1},\ 1\leq i \leq n-1,\ \alpha_n = 2\epsilon_n
$$
Il s'en suit:
$$
(n_i)_{0\leq i\leq n} = (1,2,2,2,\ldots, 2, 1),
$$
$$
(n_i^\vee)_{0\leq i\leq n} = (2,2,2,2,\ldots, 2, 2),
$$
$$
(n_i^{\vee\vee})_{0\leq i\leq n} = (4,2,2,2,\ldots, 2, 4),
$$
$$
h^\vee = 2n + 2
$$
$$
k'(C_n) = 2^{2n},\ k''(C_n) = 2^{2n+2}
$$
Donc
$$
n_i^\vee k'(C_n)^{-1/h} = 1
$$
pour tout $0\leq i \leq n$. 

D'un autre c\^ot\'e, pour $1 \leq i \leq n$, 
$$
\gamma'(C_n,\alpha_i) = 
\gamma\big(\frac{n-i}{2n}\big)\gamma\big(\frac{n-i}{2n}\big)^{-1}\gamma\big(\frac{2n-2i+1}{2n}
\big)^{-1}\gamma\big(\frac{2n-2i+1}{2n}\big)
\times
$$
$$
\times
\gamma\big(\frac{2n-2i-1}{2n}\big)^{-1}\gamma\big(\frac{2n-2i+1}{2n}\big)\gamma\big(\frac{1}{2n}\big)^{2}\gamma\big(\frac{1}{2n}\big)^{-2}\gamma\big(\frac{i}{2n}\big)^{-1}\gamma\big(\frac{2n-i}{2n}\big)^{-1} =1, 
$$
en pr\'enant compte que
$$\gamma\big(\frac{i}{2n}\big)^{-1}\gamma\big(\frac{2n-i}{2n}\big)^{-1}=\gamma\big(\frac{i}{2n}\big)^{-1}\gamma\big(\frac{i}{2n}\big).$$
\\
Donc la formule $(F')$ est v\'erifi\'ee pour $C_n$. 

\bigskip

{\it Formule $(F'')$}

\bigskip

V\'erifions la formule $(F'')$ pour $1 \leqslant  i \leqslant n-1$.
$$
\gamma''(C_n,\alpha_i) = 
\gamma\big(\frac{2n-2i+2}{2(n+1)}\big)^{-1}\gamma\big(\frac{2n-2i}{2(n+1)}\big)\gamma\big(\frac{n-i}{2(n+1)}\big)^{-1} \times
$$
$$
\times
\gamma\big(\frac{i}{2(n+1)}\big)^{-1}\gamma\big(\frac{2n-i+1}{2(n+1)}\big)^{-1}\gamma\big(\frac{n-i+1}{2(n+1)}\big)$$
Sachant que : 
$$\Gamma\big(\frac{2i}{2n+2}\big)=\Gamma\big(\frac{i}{2n+2}\big)\Gamma\big(\frac{n+1+i}{2n+2}\big)\pi^{-1/2}2^{\frac{2i-2n-2}{2n+2}}$$
$$\Gamma\big(\frac{2(n-i+1)}{2n+2}\big)=\Gamma\big(\frac{n-i+1}{2n+2}\big)\Gamma\big(\frac{2n-i+2}{2n+2}\big)\pi^{-1/2}2^{\frac{-2i}{2n+2}}$$
$$\Gamma\big(\frac{2n-2i}{2n+2}\big)=\Gamma\big(\frac{n-i}{2n+2}\big)\Gamma\big(\frac{2n+1-i}{2n+2}\big)\pi^{-1/2}2^{\frac{-2i-2}{2n+2}}$$
$$\Gamma\big(\frac{2i+2}{2n+2}\big)=\Gamma\big(\frac{i+1}{2n+2}\big)\Gamma\big(\frac{l+i+2}{2n+2}\big)\pi^{-1/2}2^{\frac{2i-2n}{2n+2}},
$$
on obtient 
$$
\gamma''(C_n,\alpha_i) = 2^{\frac{-2}{n+1}} = n_i^{\vee\vee} k''(C_n)^{-1/(2n+2)},
$$
$1 \leqslant  i \leqslant n-1$. 

Ensuite, 
$$
\gamma''(C_n,\alpha_n) = 
\gamma\big(\frac{1}{n+1}\big)^{-1}\gamma\big(\frac{l}{2(n+1)}\big)^{-1}\gamma\big(\frac{1}{2(n+1)}\big)$$
Sachant que: 
$$\Gamma\big(\frac{2n}{2n+2}\big)=\Gamma\big(\frac{l}{2n+2}\big)\Gamma\big(\frac{2n+1}{2n+2}\big)\pi^{-1/2}2^{\frac{-2}{2n+2}}$$
et
$$\Gamma\big(\frac{2}{2n+2}\big)=\Gamma\big(\frac{1}{2n+2}\big)\Gamma\big(\frac{n+2}{2n+2}\big)\pi^{-1/2}2^{\frac{-2n}{2n+2}},
$$
il vient  
$$
\gamma''(C_n,\alpha_n) = 2^{\frac{n-1}{n+1}} = n_n^{\vee\vee} k''(C_n)^{-1/(2n+2)} 
$$
Donc la formule $(F'')$ est v\'erifi\'ee pour $C_n$.

\bigskip

{\it Syst\`eme de type $F_4$}
 
\bigskip

Ce syst\`eme est obtenu par "pliure" du syst\`eme $E_6$ (le graphe de Dynkin affine dual de $F_4^{(1)}$ est 
$E_6^{(2)}$). 

$h = 12$
$$
\theta = 2\alpha_1 + 3\alpha_{2}+4\alpha_{3}+ 2\alpha_{4}
$$
$$
\rho = \frac{1}{2}(11\epsilon_1 + 5\epsilon_2 + 3\epsilon_3 + \epsilon_4)
$$
Racines simples: 
$$
\alpha_0 = - \epsilon_1 - \epsilon_2;\ \alpha_i = \epsilon_i - \epsilon_{i+1},\ i = 1,2;\ 
\alpha_3 = \epsilon_4;\ \alpha_4 = \frac{1}{2}(\epsilon_1 - \epsilon_2 - \epsilon_3 - \epsilon_4)
$$
Il s'en suit:
$$
(n_i)_{0\leq i\leq 4} = (1, 2, 3, 4, 2),
$$
$$
(n_i^\vee)_{0\leq i\leq 4} = (1, 2, 3, 2, 1),
$$
$$
(n_i^{\vee\vee})_{0\leq i\leq 4} = (1, 2, 3, 1, 1/2)
$$
Donc $h^\vee = 9$, 
$$
k'(F_4) = 2^63^3 = k(E_6)
$$
$$
k''(F_4) = 2\cdot 3^3
$$

$$
\gamma'(F_4,\alpha_1) = \frac{\gamma(3/12)}{\gamma(4/12)\gamma(5/12)} = 
\gamma(E_6,\alpha_2) = 2^{1/2}3^{-1/4} = 
n_1^\vee k'(F_4)^{-1/12}
$$
$$
\gamma'(F_4,\alpha_2) = \frac{\gamma(1/12)\gamma(4/12)\gamma(5/12)^2}{\gamma(3/12)^3} = 
\gamma(E_6,\alpha_4) = 2^{-1/2}3^{3/4} = 
n_2^\vee k'(F_4)^{-1/12}
$$
$$
\gamma'(F_4,\alpha_3) = \frac{\gamma(3/12)}{\gamma(4/12)\gamma(5/12)} = 
\gamma(E_6,\alpha_2) = 2^{1/2}3^{-1/4} = 
n_3^\vee k'(F_4)^{-1/12}
$$
$$
\gamma'(F_4,\alpha_4) = \frac{\gamma(3/12)}{\gamma(1/12)\gamma(8/12)} = 
\gamma(E_6,\alpha_1) = 2^{-1/2}3^{-1/4} = 
n_4^\vee k'(F_4)^{-1/12}
$$
La formule $(F')$ pour $F_4$ est v\'erifi\'ee. 

\bigskip

{\it Formule $(F'')$}

\bigskip

On remarque que les nombres $(\alpha|\rho)$ appartiennent \`a $\frac{1}{2}\BZ$, mais pas forcement \`a 
$\BZ$, d'o\`u l'apparence de $18 = 2h^\vee$ dans certains d\'enominateurs. 

$$
\gamma''(F_4,\alpha_1) = \gamma\big(\frac{3}{18}\big)\gamma\big(\frac{4}{9}\big)^{-1}\gamma\big(\frac{7}{18}\big)^{-1}\gamma\big(\frac{3}{9}\big)^{-1}$$
Or en utilisant la "table de 3":
$$\Gamma(\frac{3}{18})=(\Gamma(\frac{1}{18})\Gamma(\frac{7}{18})\Gamma(\frac{13}{18})(2\pi)^{-1}3^{3/18-1/2}.$$
$$\Gamma(\frac{15}{18})= (\Gamma(\frac{5}{18})\Gamma(\frac{11}{18})\Gamma(\frac{17}{18})(2\pi)^{-1}3^{15/18-1/2}.$$
$$\Gamma(\frac{6}{18})= (\Gamma(\frac{2}{18})\Gamma(\frac{8}{18})\Gamma(\frac{14}{18})(2\pi)^{-1}3^{6/18-1/2}.$$
$$\Gamma(\frac{12}{18})= (\Gamma(\frac{4}{18})\Gamma(\frac{10}{18})\Gamma(\frac{16}{18})(2\pi)^{-1}3^{12/18-1/2}.$$
et "la table de 2" :
$$\Gamma(\frac{2}{18})= \pi^{-1/2}2^{1/9-1} \Gamma(\frac{1}{18})\Gamma(\frac{10}{18})$$
$$\Gamma(\frac{4}{18})= \pi^{-1/2}2^{2/9-1} \Gamma(\frac{2}{30})\Gamma(\frac{17}{30})$$
$$\Gamma(\frac{8}{18})= \pi^{-1/2}2^{4/9-1} \Gamma(\frac{4}{30})\Gamma(\frac{19}{30})$$
$$\Gamma(\frac{10}{18})= \pi^{-1/2}2^{5/9-1} \Gamma(\frac{7}{30})\Gamma(\frac{22}{30})$$
$$\Gamma(\frac{14}{18})= \pi^{-1/2}2^{7/9-1} \Gamma(\frac{8}{30})\Gamma(\frac{23}{30})$$
$$\Gamma(\frac{16}{18})= \pi^{-1/2}2^{8/9-1} \Gamma(\frac{11}{30})\Gamma(\frac{26}{30}),
$$
on obtient 
$$
\gamma''(F_{4}, \alpha_{1})=3^{-1/3}2^{8/9} = n_1^{\vee\vee}k''(F_4), 
$$
la formule $(F'')$ est donc v\'erifi\'ee pour $i=1$. 

Ensuite, 
$$
\gamma''(F_{4}, \alpha_{2})=\gamma\big(\frac{3}{18}\big)^{-1}
\gamma\big(\frac{2}{9}\big)^{-1}\gamma\big(\frac{5}{9}\big)^{-1}
\gamma\big(\frac{1}{18}\big)\gamma\big(\frac{7}{18}\big)
$$
En utilisant les "tables de 3 et de 2", on obtient 
$$
\gamma''(F_{4}, \alpha_{2})=2^{-1/9}3^{2/3} = n_2^{\vee\vee}k''(F_4). 
$$
la formule $(F'')$ est donc v\'erifi\'ee pour $i=2$. 

Ensuite,
$$
\gamma''(F_{4}, \alpha_{3})=\gamma\big(\frac{1}{18}\big)^{-1}\gamma\big(\frac{7}{18}\big)^{-1}
\gamma\big(\frac{5}{18}\big)^{-1}\gamma\big(\frac{1}{3}\big)\gamma\big(\frac{5}{9}\big)\gamma\big(\frac{1}{9}\big)\gamma\big(\frac{2}{9}\big).
$$

Or en utilisant les "tables de 3 et de 2", on obtient  
$$\gamma''(F_{4}, \alpha_{3})=2^{-1/9}3^{-1/3} = n_3^{\vee\vee}k''(F_4), 
$$
ce qui v\'erifie $(F'')$ pour $i=3$. 

Enfin,
$$
\gamma''(F_{4}, \alpha_{4})=\gamma\big(\frac{1}{18}\big)^{-1}\gamma\big(\frac{2}{18}\big)^{-1}
\gamma\big(\frac{6}{18}\big)^{-1}\gamma\big(\frac{10}{18}\big)^{-1}\gamma\big(\frac{11}{18}\big)^{-1}
\gamma\big(\frac{3}{18}\big)\gamma\big(\frac{5}{18}\big)\gamma\big(\frac{4}{18}\big).
$$
Or en utilisant les "tables de 3 et de 2", on en d\'eduit 
$$
\gamma''(F_{4}, \alpha_{4})=2^{-10/9}3^{-1/3} = n_4^{\vee\vee}k''(F_4), 
$$
ce qui v\'erifie $(F'')$ pour $i=4$ et ach\`eve la v\'erification de $(F'')$ pour le syst\`eme $F_4$.  

\bigskip

{\it Syst\`eme de type $G_2$}
 
\bigskip

Ce syst\`eme est obtenu par "pliure" du syst\`eme $D_4$ (le graphe de Dynkin affine dual de $G_2^{(1)}$ est 
$D_4^{(3)}$).

$h = 6$
$$
\theta = 3\alpha_1 + 2\alpha_2
$$
Racines simples:
$$
\alpha_{0}= - \epsilon_{1} - \epsilon_{2}+2\epsilon_{3},\ 
\alpha_{1}= \epsilon_{1}-\epsilon_{2},\ \alpha_{2}= -2\epsilon_{1}+\epsilon_{2}+\epsilon_{3}.
$$
Donc
$$
(n_0, n_1, n_2) = (1, 3, 2), 
$$
$$
(n_0^\vee, n_1^\vee, n_2^\vee) = (3, 3, 6), 
$$
$$
(n_0^{\vee\vee}, n_1^{\vee^\vee}, n_2^{\vee\vee}) = (9, 3, 18). 
$$

$h^\vee = 12$. 

Il s'en suit:
$$
k'(G_2) = 2^23^6
$$
$$
k''(G_2) = 2^63^{21}
$$

$$
\gamma'(G_2,\alpha_1) = \frac{\gamma(2/6)}{\gamma(1/6)\gamma(4/6)} = 
\gamma(D_4,\alpha_1) = 2^{-1/3} = 
n_1^\vee k'(G_2)^{-1/6}
$$
$$
\gamma'(G_2,\alpha_2) = \frac{\gamma(1/6)^2\gamma(4/6)}{\gamma(2/6)^3} = 
\gamma(D_4,\alpha_2) = 2^{2/3} = 
n_2^\vee k'(G_2)^{-1/6}
$$
Ceci prouve $(F')$ pour $G_2$. 

\bigskip

V\'erifions maintenant la formule (F''). 

On a 
$$
\gamma''(G_2,\alpha_1) = \gamma\big(\frac{1}{12}\big)^{-2}\gamma\big(\frac{3}{12})^{3}\gamma\big(\frac{4}{12}\big)
\gamma\big(\frac{5}{12}\big)^{-1}\gamma\big(\frac{6}{12})^{-3}
$$
En utilisant: 
$$ 
2\pi 3^{1/4}\Gamma(\frac{3}{12})=\Gamma(\frac{1}{12})\Gamma(\frac{5}{12})\Gamma(\frac{9}{12}), 
$$
puis: 
 $$2^{1/6}\pi^{1/2}\Gamma(\frac{10}{12})=\Gamma(\frac{5}{12})\Gamma(\frac{11}{12})$$
 $$2^{1/3}\pi^{1/2}\Gamma(\frac{8}{12})=\Gamma(\frac{4}{12})\Gamma(\frac{10}{12})$$
et 
$$
2^{5/6}\pi^{1/2}\Gamma(\frac{2}{12})=\Gamma(\frac{1}{12})\Gamma(\frac{7}{12})
$$
puis :
$$
\Gamma(\frac{2}{12})\Gamma(\frac{10}{12})=2\pi, 
$$
on a bien: 
$$
\gamma''(G_2,\alpha_1) = \gamma\big(\frac{1}{12}\big)^{-2}\gamma\big(\frac{3}{12})^{3}\gamma\big(\frac{4}{12}\big)
\gamma\big(\frac{5}{12}\big)^{-1}\gamma\big(\frac{6}{12})^{-3}=2^{-1/2}3^{-3/4} = n_1^\vee k''(G_2)^{-1/12}
$$
La formule $(F'')$ est v\'erifi\'ee pour $i=1$. 

Ensuite, 
$$
\gamma''(G_2,\alpha_2) = 
\gamma\big(\frac{1}{12}\big)\gamma\big(\frac{3}{12})^{-1}\gamma\big(\frac{6}{12}\big)
\gamma\big(\frac{4}{12}\big)^{-1}= \gamma\big(\frac{1}{12}\big)\gamma\big(\frac{3}{12})^{-1}\gamma\big(\frac{4}{12}\big)^{-1}
$$
En utilisant
$$ 
2\pi 3^{1/4}\Gamma(\frac{3}{12})=\Gamma(\frac{1}{12})\Gamma(\frac{5}{12})\Gamma(\frac{9}{12}),
$$
puis 
$$
2^{1/6}\pi^{1/2}\Gamma(\frac{10}{12})=\Gamma(\frac{5}{12})\Gamma(\frac{11}{12})
$$
$$2^{1/3}\pi^{1/2}\Gamma(\frac{8}{12})=\Gamma(\frac{4}{12})\Gamma(\frac{10}{12}),
$$
on a: 
$$
\gamma''(G_2,\alpha_2) = 
\gamma\big(\frac{1}{12}\big)\gamma\big(\frac{3}{12})^{-1}\gamma\big(\frac{4}{12}\big)^{-1}=2^{1/2}3^{1/4} =  
n_2^\vee k''(G_2)^{-1/12}
$$
La formule $(F'')$ est v\'erifi\'ee pour $i=2$. 

Donc la formule $(F'')$ est v\'erifi\'ee pour le syst\`eme $G_{2}$.

Ainsi les formules $(F')$ et $(F'')$ sont v\'erifi\'ees pour les syst\`emes du type B, C, F, G, ce qui ach\`eve la d\'emonstration 
du Th\'eor\`eme 1.2.

\bigskip\bigskip



\centerline{\bf Bibliographie}

\bigskip\bigskip

[ABFKR] C.Ahn, P.Baseilhac, V.A.Fateev, C.Kim, C.Rim, Reflection amplitudes in non-simply laced 
Toda theories and thermodynamic Bethe Ansatz, {\it Phys. Let.} {\bf B481} (2000), 114 - 124. 

[B] N.Bourbaki, Groupes et alg\`ebres de Lie, Chapitres. IV - VI. 

[CS] V.Cohen-Aptel, V.Schechtman, Produits Gamma et vecteurs propres des matrices de Cartan, 
arXiv:1010.5945. 

[F] V.A.Fateev, Normalization factors, reflection amplitudes and integrable systems, hep-th/0103014. 

[K] V.G.Kac, Infinite dimensional Lie algebras. 

\bigskip\bigskip

Institut de Math\'ematiques, Universit\'e Paul Sabatier, 118 Route de Narbonne, 31062 Toulouse, France

vero.aptel\verb+@+free.fr

\end{document}